# Zeros of the Alternating Zeta Function on the Line $\Re(s)=1$

Jonathan Sondow

The Riemann zeta function $\zeta(s)$ is defined as the analytic continuation of the Dirichlet series

$$\zeta(s) = \sum_{n=1}^{\infty} \frac{1}{n^s},$$

which converges if $\Re(s) > 1$. The zeta function is holomorphic in the complex plane ($= s$-plane), except for a simple pole at $s = 1$. The alternating zeta function $\zeta^*(s)$ is defined as the analytic continuation of the Dirichlet series

$$\zeta^*(s) = \sum_{n=1}^{\infty} \frac{(-1)^{n-1}}{n^s},$$

which converges if $\Re(s) > 0$. The two functions are related for all complex $s$ by the identity

(1) $$\zeta^*(s) = (1 - 2^{1-s})\zeta(s),$$

which is easily established, first for $\Re(s) > 1$ by combining terms in the convergent Dirichlet series, and then by using analytic continuation to extend the result to the entire $s$-plane (see Theorem 13.11 in [1]). The factor $(1 - 2^{1-s})$ has a simple zero at $s = 1$ that cancels the pole of $\zeta(s)$, so (1) shows that $\zeta^*(s)$ is an entire function of $s$ and that it vanishes at each zero of the factor $(1 - 2^{1-s})$ with the exception of $s = 1$, at which point $\zeta^*(1) = \log 2$. The other zeros of the factor $(1 - 2^{1-s})$ lie on the line $\Re(s) = 1$, and occur at the points $s = 1 + 2k\pi i/\log 2$, where $k$ is a nonzero integer.

This note deduces the vanishing of $\zeta^*(s)$ at the zeros of $(1 - 2^{1-s})$, except for $s = 1$, without using (1). Instead, we use identity (4) relating the partial sums

$$\zeta_N(s) = \sum_{n=1}^{N} \frac{1}{n^s}, \quad \zeta_N^*(s) = \sum_{n=1}^{N} \frac{(-1)^{n-1}}{n^s},$$

which are entire functions of $s$. To obtain (4), take an even number of terms in each of these finite sums and subtract. The terms with odd index $n$ cancel, and we find that

(2) $$\zeta_{2N}^*(s) - \zeta_{2N}(s) = -2^{1-s}\zeta_N(s).$$

But

$$\zeta_N(s) = \zeta_{2N}(s) - \sum_{n=N+1}^{2N} n^{-s}.$$

Using this in (2) and rearranging terms, we obtain the identity

(3) $$\zeta_{2N}^*(s) = (1 - 2^{1-s})\zeta_{2N}(s) + 2^{1-s} \sum_{n=N+1}^{2N} n^{-s},$$

which holds for all complex $s$.

The last sum in (3) can be written as

$$\sum_{n=N+1}^{2N} n^{-s} = \sum_{k=1}^{N} (N+k)^{-s} = N^{1-s} \frac{1}{N} \sum_{k=1}^{N} (1 + \tfrac{k}{N})^{-s},$$

and we recognize the factor multiplying $N^{1-s}$ as a Riemann sum for the integral

$$\int_0^1 (1+x)^{-s} dx = \begin{cases} \frac{1 - 2^{1-s}}{s-1} & \text{if } s \neq 1, \\ \log 2 & \text{if } s = 1. \end{cases}$$

Therefore, we introduce the difference

$$d_N(s) = \int_0^1 (1+x)^{-s} dx - \frac{1}{N} \sum_{k=1}^{N} (1 + \tfrac{k}{N})^{-s}$$

and rewrite (3) as follows:

(4) $$\zeta_{2N}^*(s) = (1 - 2^{1-s})\zeta_{2N}(s) + (2N)^{1-s} \left\{ \int_0^1 (1+x)^{-s} dx - d_N(s) \right\}.$$

Note that $d_N(s)$ tends to 0 as $N \to \infty$. For example, when $s = 1$ equation (4) becomes

$$\zeta^*_{2N}(1) = \log 2 - d_N(1)$$

and when $N \to \infty$ it gives the classic result $\zeta^*(1) = \log 2$ for the alternating harmonic series. Equation (4) is valid for all $s$, and we use it to show that $\zeta^*(s)$ vanishes at each zero of $1 - 2^{1-s}$ other than $s = 1$. In fact, if we take $s = 1 + it$, where $t$ is not zero, then (4) becomes

$$\zeta^*_{2N}(1+it) = (1 - 2^{-it})\zeta_{2N}(1+it) + (2N)^{-it}\left\{\tfrac{1-2^{-it}}{it} - d_N(1+it)\right\}.$$

If $t$ is chosen so that $2^{-it} = 1$, that is, if $t = 2k\pi / \log 2$ for any integer $k \neq 0$, we get

$$\zeta^*_{2N}(1+it) = -N^{-it} d_N(1+it).$$

The factor $N^{-it}$ has modulus 1 and $d_N(1+it) \to 0$ as $N \to \infty$, hence $\zeta^*(1+it) = 0$. ●

**Note.** In [2] we show that the Riemann hypothesis is related to the rate at which $d_N(s)$ converges to zero in the critical strip $0 < \Re(s) < 1$.

**ACKNOWLEDGEMENT.** The author thanks the referee for suggesting the changes that led to the present version of the paper.

*209 West 97th St., New York, NY 10025*
*jsondow@alumni.princeton.edu*